\documentclass[12pt]{article} 
\pdfoutput=1
\usepackage[authoryear]{natbib}
\usepackage{graphicx,epstopdf}
\usepackage{setspace}
\usepackage[margin=2.5cm]{geometry}
\usepackage{array}
\usepackage{hyperref}

\usepackage{pdflscape}
\usepackage{graphicx}
\usepackage[utf8]{inputenc}
\usepackage[english]{babel}
\usepackage{amsmath, amsthm, amssymb, amsfonts}

\usepackage{ifpdf}
\ifpdf
  \DeclareGraphicsExtensions{.pdf,.png,.jpg}
\else
  \DeclareGraphicsExtensions{.eps}
\fi

\usepackage{natbib}
\setcitestyle{authoryear,open={(},close={)}}

\usepackage{CJKutf8}

\usepackage{framed}

\usepackage{arydshln}

\usepackage{enumitem}

\usepackage{float}

\usepackage{hyperref}

\usepackage{comment}

\usepackage{tikz}
\usetikzlibrary{arrows}

\usepackage{multirow}

\usepackage{amsmath, amssymb, amsthm, amsfonts, bbm, mathtools, dsfont}
\DeclareMathOperator{\vech}{vech}

\newcommand{\EE}{\mathbb{E}}

\newcommand{\NN}{\mathbb{N}}

\newcommand{\RR}{\mathbb{R}}

\newcommand{\ZZ}{\mathbb{Z}}

\newcommand{\cC}{\mathcal{C}}

\newcommand{\cI}{\mathcal{I}}

\newcommand{\cN}{\mathcal{N}}

\newcommand{\cT}{\mathcal{T}}

\newcommand{\bG}{\mathbf{G}}

\newcommand{\bI}{\mathbf{I}}

\newcommand{\bK}{\mathbf{K}}

\newcommand{\bT}{\mathbf{T}}
\newcommand{\bU}{\mathbf{U}}

\newcommand{\bV}{\mathbf{V}}

\newcommand{\bX}{\mathbf{X}}
\newcommand{\bY}{\mathbf{Y}}
\newcommand{\bZ}{\mathbf{Z}}

\newcommand{\by}{\mathbf{y}}
\newcommand{\bv}{\mathbf{v}}
\newcommand{\bg}{\mathbf{g}}

\DeclarePairedDelimiter{\floor}{\lfloor}{\rfloor}
\DeclarePairedDelimiter{\ceil}{\lceil}{\rceil}

\newcommand{\indicator}{\mathds{1}}
\DeclareMathOperator{\var}{Var}

\newcommand{\Normal}{\text{N}}

\newcommand{\simiid}{\mathrel{\stackrel{\makebox[0pt]{\mbox{\normalfont\tiny iid}}}{\sim}}}

\newcommand\numberthis{\addtocounter{equation}{1}\tag{\theequation}}

\newtheorem{theorem}{Theorem}[section]

\newtheorem{lemma}[theorem]{Lemma}

\theoremstyle{definition}
\newtheorem{definition}{Definition}[section]
\newtheorem{assumption}{Assumption}[section]

\usepackage[stable]{footmisc}

\title{Minimax Risk and Uniform Convergence Rates for Nonparametric Dyadic Regression}

\begin{small}
\author{Bryan S. Graham\footnote{Department of Economics, University of California - Berkeley and the National Bureau of Economic Research, e-mail: \texttt{bgraham@econ.berkeley.edu}}, Fengshi Niu\footnote{Department of Economics, University of California - Berkeley,  e-mail: \texttt{fniu@berkeley.edu}}, James L. Powell\footnote{Department of Economics, University of California - Berkeley, e-mail: \texttt{powell@econ.berkeley.edu}}, \thanks{We thank seminar audiences at University of California - Berkeley for helpful feedback. We also thank Matias Cattaneo, Michael Jansson and Harold Chiang for useful comments and discussion. All the usual disclaimers apply. Financial support from the National Science Foundation (SES \#1357499, SES \#1851647) is gratefully acknowledged.}} 
\end{small}

\date{\textsc{Initial Draft: July 2019, This Draft: March 2021}}

\begin{document}
\maketitle

\singlespacing
\begin{abstract}
\noindent \singlespacing
Let $i=1,\ldots,N$ index a simple random sample of units drawn from some large population. For each unit we observe the vector of regressors $X_{i}$ and, for each of the $N\left(N-1\right)$ ordered pairs of units, an outcome $Y_{ij}$. The outcomes $Y_{ij}$ and $Y_{kl}$ are independent if their indices are disjoint, but dependent otherwise (i.e., ``dyadically dependent''). Let $W_{ij}=\left(X_{i}',X_{j}'\right)'$; using the sampled data we seek to construct a nonparametric estimate of the mean regression function $g\left(W_{ij}\right)\overset{def}{\equiv}\mathbb{E}\left[\left.Y_{ij}\right|X_{i},X_{j}\right].$    

We present two sets of results. First, we calculate lower bounds on the minimax risk for estimating the regression function at (i) a point and (ii) under the infinity norm. Second, we calculate (i) pointwise and (ii) uniform convergence rates for the dyadic analog of the familiar Nadaraya-Watson (NW) kernel regression estimator. We show that the NW kernel regression estimator achieves the optimal rates suggested by our risk bounds when an appropriate bandwidth sequence is chosen. This optimal rate differs from the one available under iid data: the effective sample size is smaller and $d_W=\mathrm{dim}(W_{ij})$ influences the rate differently.
\end{abstract}



\newpage
\setcounter{page}{1}

\onehalfspacing
\section{Introduction}
Let $i=1,\ldots,N$ index a simple random sample of units drawn from some large population. For each unit we observe the vector of regressors $X_{i}$ and, for each of the $N\left(N-1\right)$ ordered pairs of units, or \emph{directed dyads}, we observe the ``dyadic” outcome $Y_{ij}$ (e.g., total exports from country $i$ to country $j$). The outcomes $Y_{ij}$ and $Y_{kl}$ are independent if their indices are disjoint, but dependent otherwise (e.g., exports from Japan to Korea may covary with those from Japan to Vietnam). 

Let $W_{ij}=\left(X_{i}',X_{j}'\right)'$; using the sampled data we seek to construct a nonparametric estimate of the mean regression function 
\begin{equation}\label{eq: mean_reg}
    g\left(W_{ij}\right)\overset{def}{\equiv}\mathbb{E}\left[\left.Y_{ij}\right|X_{i},X_{j}\right].    
\end{equation}

We present two sets of results. First, we calculate lower bounds on the minimax risk for estimating the regression function at (i) a point and (ii) under the infinity norm. Second, we calculate (i) pointwise and (ii) uniform convergence rates for the dyadic analog of the familiar Nadaraya-Watson (NW) kernel regression estimator. We show that the NW kernel regression estimator achieves the optimal rates suggested by our risk bounds when an appropriate bandwidth sequence is chosen.

Analogous results are widely available in the i.i.d. setting. For nonparametric regression risk bounds see, for example,  \cite{stone1980,stone1982} and \cite{ibragimov1982, ibragimov1984}. \cite{tsybakov2008} provides a masterful synthesis of these results, from which we draw in formulating our own proofs.

Uniform convergence of kernel averages with i.i.d. data, as well as stationary strong mixing data, have been studied by, for example, \cite{newey1994} and \cite{hansen2008} respectively. The latter paper includes additional references to the extensive literature in this area. Our uniform convergence proofs build upon those of \cite{hansen2008}. Nonparametric density estimation with dyadic data was first considered by \cite{graham2019}; \cite{chiang2019} present uniform convergence results for dyadic density estimators.\footnote{It is possible that the methods of inference presented in \cite{chiang2019} could be adapted to our setting.}

Our results provide insight in the structure of dyadic nonparametric estimation problems. Our minimax risk bounds suggest that, $N$, the number of units, \emph{not} $n \overset{def}{\equiv} N \times (N-1)$, the number of dyadic outcomes, is the relevant ``sample size'' for dyadic estimation problems. This is consistent with the long standing intuition among empirical researchers that dyadic dependence makes inference less precise (see \cite{aronow_samii_assenova_2017} and the references cited therein), as well as with a small, but growing, number of more formal rates-of-convergence results \citep[cf.,][]{graham2020network}.

More surprisingly, we find that the relevant dimension of our estimation problem is just $d_X = \mathrm{dim}(X_i)$, not $d_W = 2d_X$. We provide two intuitions for this fact. The first, described below, stems from the thought experiment underlying our minimax risk bound calculations. The second, arises from the fact that the H\'ajek projection of the NW estimator has a ``partial-mean-like" structure. As is well known, averaging over the marginal distribution of some regressors, while holding the remaining ones fixed, improves rates-of-convergence \citep[e.g.,][]{newey1994, linton1995}.

\cite{graham2020network} surveys empirical studies in economics utilizing dyadic data. Interest in, as well as the availability of, such data are growing in economics, other academic fields, and in enterprise settings. This paper provides an initial set of results for nonparametric regression with dyadic data. These results are, of course, of direct interest. They should, as has been true with their i.i.d. predecessors, also be useful for proving consistency of two-step semiparametric M-estimators under dyadic dependence (see \cite{chiang2019} for some results on double machine learning with dyadic data).

\section{Lower Bounds on the Minimax Risk} \label{sec: risk_bounds}
Let $i=1,\ldots,N$ index a simple random sample of units drawn from some large population. The econometrician observes the vector of regressors, $X_{i}$, for each sampled unit as well as the scalar outcome, $Y_{ij}$, for each directed pair of sampled units (i.e., each directed dyad). Let $\bZ_N = (X_1, \ldots, X_N, Y_{ij}, 1 \leq i\neq j \leq N)$ be the observable data when $N$ units are sampled. The regression function of interest is \eqref{eq: mean_reg} above. The goal is to construct a nonparametric estimate of $g:\mathbb{R}^{d_{W}}\rightarrow\mathbb{R}$ where $d_{W}=2d_{x}$.

We assume that $Y_{ij}$ is generated according to the following conditionally independent dyad (CID) model \cite[cf.,][Section 3.3]{graham2020network}.
\begin{equation}
    Y_{ij}=h(X_{i},X_{j},U_{i},U_{j},V_{ij}). \label{eq: graphon}
\end{equation}
Random sampling ensures that $\left(X_{i},U_{i}\right)$ is independent and identically distributed for $i=1,\ldots,N$. We further assume that $\left\{ \left(V_{ij},V_{ji}\right)\right\} _{1\leq i<j\leq N}$ are i.i.d. and indepenent of $\mathbf{X}=\left(X_{1},\ldots,X_{N}\right)'$ and $\mathbf{U}=\left(U_{1},\ldots,U_{N}\right)$. Here $h$ is an unknown function, often called the \emph{graphon}.
This set-up, which can also be derived as an implication of more primitive exchangeability assumptions, has the following implications (see \cite{graham2020network,graham2020sparse} for additional discussion):
\begin{enumerate}
\item The $Y_{ij}$ are relatively exchangeable given the $W_{ij}$. Namely, the conditional distribution of $\bY$ is invariant across permutations of the indices $\sigma : \NN \rightarrow \NN$ satisfying the restriction $[W_{\sigma(i)\sigma(j)}] \overset{d}{=} [W_{ij}]$:
\[
[Y_{ij}] \overset{d}{=} [Y_{\sigma(i)\sigma(j)}].
\]
\item $Y_{ij}$ and $Y_{kl}$ are independent if their indices are disjoint.
\item $Y_{ij}$ and $Y_{kl}$ are dependent (unconditionally or conditionally given $X_1,\ldots,X_N$) if they share at least one index in common.
\end{enumerate}

The statistical problem is to estimate the regression function $g$ when the only prior restriction on it is that it belongs to the H\"older class of functions.
\begin{definition}
\textsc{(H\"older Class)} Given a vector $s = (s_1, \ldots, s_{d})$, define $|s| = s_1 + \cdots + s_{d}$ and \[D^s = \frac{\partial^{s_1 + \cdots + s_{d}}}{\partial^{s_1} w_1 \cdots \partial^{s_{d}} w_{d}}.\]
Let $\beta$ and $L$ be two positive numbers. The \emph{H\"older class} $\Sigma(\beta, L)$ on $\RR^{d}$ is defined as the set of $l = \floor{\beta}$ times differentiable functions $g:\RR^{d} \rightarrow \RR$ whose partial derivative $D^s g$ satisfies 
\[
|D^s g(w) - D^s g(w')| \leq L||w - w'||_{\infty}^{\beta - l}, \quad \forall w, w'\in \RR^{d}
\]
for all $s$ such that $|s| = \floor{\beta}$. $\floor{\beta}$ denotes the greatest integer strictly less than the real number $\beta$.
\end{definition}

An estimator $\hat{g}_N$ is a function $w \mapsto \hat{g}_N(w) = \hat{g}_N(w, \bZ_N)$ measurable with respect to $\bZ$. Our first result establishes a lower bound on the minimax risk for estimating the regression function at a single point and under the infinity norm. We state this result under a Gaussian error assumption, which simplifies the proof. 

\begin{theorem}
\label{thm:Minimax}
\textsc{(Minimax Risk Lower Bound)} Suppose that $\beta > 0$ and $L > 0$; $X_i$ is continuously distributed on $\RR^{d_X}$ with density $f$ and $\sup_x f(x) \leq B_3 < \infty$; and $Y_{ij}$ is generated according to the following nonparametric regression model:
\begin{align*}
Y_{ij} & = g\left(W_{ij}\right) + e_{ij}, \quad i \neq j,
\end{align*}
with $e_{ij} = U_i + U_j + V_{ij}$, $U_i \simiid \Normal(0, 1)$, and $V_{ij}\simiid \Normal(0, 1)$, then 
\begin{enumerate}[label=(\roman*)]
\item 
For all $w \in \RR^{d_W}$, 
\begin{align*}
\liminf_{N\rightarrow \infty} \inf_{\hat{g}_N} \sup_{g \in \Sigma(\beta, L)}
\EE_{g} \left[N^{\frac{2\beta}{2\beta + d_X}}\left(\hat{g}_N(w) - g(w)\right)^2\right] \geq c_1,
\end{align*}
where $c_1 > 0$ depends only on $\beta$ and $L$.

\item
\begin{align*}
\liminf_{N\rightarrow \infty} \inf_{\hat{g}_N} \sup_{g \in \Sigma(\beta, L)}
\EE_{g} \left[\left(\frac{N}{\ln N}\right)^{\frac{2\beta}{2\beta + d_X}}\left|\left|\hat{g}_N - g\right|\right|_{\infty}^2\right] \geq c_2,
\end{align*}
where $c_2 > 0$ also depends only on $\beta$ and $L$.
\end{enumerate}
\end{theorem}

Our proof follows the general recipe outlined in Chapter 2 of \cite{tsybakov2008}. The lower bound at a point is based on Le Cam's method of two hypotheses. The lower bound under the infinity norm is based on Fano's method of multiple hypotheses. 

The key, and novel, step in our proof involves constructing hypotheses close enough to one other in terms of Kullback-Leibler (KL) divergence while being at the same time different enough in terms of the target regression function.

An essential feature of our construction is additive separability of the regression functions. In the hypotheses we consider, $Y_{ij} = k(X_i) + k(X_j) + U_i + U_j + V_{ij}$. Next suppose we also observe $T_i \overset{def}{\equiv} k(X_i) + U_i$. Observe that $(X_i, T_i, i = 1,\ldots, N)$ is sufficient with respect to $(X_i, T_i, i = 1,\ldots, N, Y_{kl}, 1\leq k\neq l \leq N)$ for the parameter $k$. 

It is well-known that the optimal rates of convergence for estimating $k$ using iid data $(X_i, T_i, i = 1,\ldots, N)$ are  $N^{-\frac{\beta}{2\beta + d_X}}$ pointwise and $\left(\frac{N}{\ln N}\right)^{-\frac{\beta}{2\beta + d_X}}$ for the infinity norm. We expect the rates for estimating $g$ to be no faster than these. The proof of Theorem \ref{thm:Minimax} makes this intuition rigorous.

Relative to its iid counterpart, there are two distinctive features of Theorem \ref{thm:Minimax}. First, the relevant sample size is not the number of observed dyadic outcomes $n = N \times (N-1)$, but instead the number of sampled units, $N$. Dependence across outcomes sharing indices in common is strong enough to slow down the feasible rate of convergence. Second, although the regression function has $d_W = 2d_X$ arguments, the relevant dimension reflected in the rate of convergence result is just $d_X$ (i.e., just half of what might naively be expected).

The form of our constructed hypotheses provides one intuition for this second finding: clearly the relevant dimension of the problem of estimating $k(x)$ is just $d_X$. Relatedly this finding is consistent with those of \cite{linton1995} in their analysis of additively separable, but otherwise nonparametric, regression functions (see also \cite{newey1994}).

The pairwise structure of dyadic data results in apparent data abundance (sample $N$ agents, but observe $O(N^2)$ outcomes!). This abundance is both illusory, in the sense that the effective sample size is indeed just $N$, and real, in the sense that availability of the pairwise outcome data allows for an effective reduction in the dimensionality of the problem via partial mean like average (as in \cite{newey1994} and \cite{linton1995}  in a different context).

\section{Kernel Estimator of Dyadic Regression} \label{sec: NW_reg}
In this section we study the properties of a specific nonparametric regression estimator. Namely, the dyadic analog of the well-known Nadaraya-Watson (NW) kernel regression estimator. While our results are specific to this estimator, they could, for example, be extended to apply to local linear regression \citep[e.g.,][]{hansen2008}.

The dyadic NW kernel regression estimator is
\begin{equation}
	\hat{g}_N(w) := \frac{\sum_{1 \leq i \neq j \leq N} K_{ij, N}(w) Y_{ij}}{\sum_{1 \leq i \neq j \leq N} K_{ij, N}(w)},
\end{equation}
where
\[
K_{ij, N}(w) := \frac{1}{h_N^{d_W}}K\left(\frac{W_{ij} - w}{h_N}\right),
\]
$K$ is a fixed multivariate kernel function, and $h_N$ is a vanishing bandwidth sequence.  

We first develop a sequence of results useful for bounding the variance of kernel objects of the form
\begin{align}
\hat{\Psi}_N(w)
& := \frac{1}{N(N-1)}\sum_{1\leq i \neq j\leq N} Y_{ij} K_{ij, N}(w) 
\end{align}
and then apply these results to the NW regression estimator. We then bound the NW estimator's bias and combine the two sets of results to formulate a risk bound.

\subsection{Variance Bound and Uniform Convergence}
Here we are interested in bounding the deviation of $\hat{\Psi}_N(w)$ from its mean. We begin with a presentation of our maintained assumptions.
\begin{assumption}[\textsc{Model}] The data generating process is as described in Section \ref{sec: risk_bounds} with
\label{a1}
\begin{enumerate}[label=(\roman*)]
\item $X_i$ continuously distributed with marginal density $f(x)$ s.t. $\sup_{x\in \RR^{d_X}}f(x) \leq B_3 < \infty$;
\item $\sup_{x_1, x_2 \in \RR^{d_X}} \EE\left[|Y_{12}|^2\big|(X_1, X_2)=(x_1, x_2)\right] \cdot f(x_1) f(x_2) \leq B_4 < \infty$, \\$\sup_{x_1, x_2, x_3 \in \RR^{d_X}} \EE\left[|Y_{12}Y_{13}|\big|(X_1, X_2, X_3)=(x_1, x_2, x_3)\right] \cdot f(x_1) f(x_2) f(x_3) \leq B_5 < \infty$.
\end{enumerate}
\end{assumption}
Condition (i) is a standard condition in the context of kernel estimation, while (ii) ensures that various second moments appearing in our variance calculations are finite.
\begin{assumption}[\textsc{Kernel, Part A}]
\label{a2}
$\sup_{w \in \RR^{d_W}} |K(w)| \leq K_{\text{max}} < \infty$, 
$\int_{w\in \RR^{d_W}} |K(w)|\mathrm{d}w \leq B_1 < \infty$, 
and $\sup_{x\in \RR^{d_X}}\int |K(x, x')|\mathrm{d}x' \leq B_2 < \infty$.
\end{assumption}
Assumption \ref{a2} is satisfied by many widely-used multivariate kernel functions. Our first result holds under Assumptions \ref{a1} and \ref{a2}.
\begin{theorem}[\textsc{Variance Bound}]
\label{thm:VarianceBound}
Under Assumptions \ref{a1} and \ref{a2}, and the bandwidth condition $Nh_N^{d_X} \rightarrow \infty$ as $N\rightarrow \infty$, there exists a constant $M_0 < \infty$ such that for $N$ sufficiently large
\[
\var\left(\hat{\Psi}_N(w)\right) \leq \frac{M_0}{Nh_N^{d_X}}
\]
for all $w\in \RR^{d_W}$.
\end{theorem}
A proof is available in the appendix. Mirroring our risk bound results, two features of Theorem \ref{thm:VarianceBound} merit comment. First, $N$ not $n = N \times (N-1)$ appears in the denominator. This is due to the effects of dependence across dyads sharing units in common. Second, the relevant dimension of the problem is $d_X$, not $d_W=2d_X$, this reflects the U-statistic like structure of kernel weighted averages and the partial mean like averaging this structure induces.

To establish uniform convergence, we need additional moment conditions on $Y_{ij}$ as well as some smoothness conditions on the kernel $K$. As in \cite{hansen2008}, we require the kernel to either have bounded support and be Lipschitz or have bounded derivatives and an integrable tail. See \cite{hansen2008} for additional discussion about these conditions. As with Assumption \ref{a2} above, most commonly used kernels satisfy these conditions.
\begin{assumption}[\textsc{Regularity Condition}]
\label{a3}
\begin{enumerate}[label=(\roman*)]
\item 
For some $s>2$, $\EE|Y_{12}|^s < \infty$ and 
$\sup_{x_1, x_2 \in \RR^{d_X}} \EE\left[|Y_{12}|^s\big|(X_1, X_2)=(x_1, x_2)\right] \cdot f(x_1, x_2) \leq B_{4,s} < \infty$;
\item For some $\Lambda_1 < \infty$ and $L < \infty$, either (a) or (b) holds
\begin{itemize}
\item[(a)] $K(w) = 0$ for $||w|| > L$, and $|K(w) - K(w')| \leq \Lambda_1||w-w'||$ for all $w, w'\in \RR^{2d}$
\item[(b)] $K(w)$ is differentiable, $\left|\left|\frac{\partial}{\partial w} K(w)\right|\right| \leq \Lambda_1$, where $\left|\left|\frac{\partial}{\partial w} K(w)\right|\right| = \left|\left|\left(\frac{\partial}{\partial w_1} K(w), \ldots \frac{\partial}{\partial w_{2d}} K(w)\right)\right|\right|_{\infty}$, and for some $\nu > 1$, $\left|\left|\frac{\partial}{\partial w} K(w)\right|\right| \leq \Lambda_1 ||w||^{-\nu}$ for $||w|| > L$.
\end{itemize}
\end{enumerate}
\end{assumption}
Part (ii) coincides with Assumption 3 in \cite{hansen2008}. This assumption implies that for all $||w_1 - w_2|| \leq \delta \leq L$,
\begin{align*}
|K(w_2) - K(w_1)| \leq \delta K^*(w_1),
\end{align*}
where $K^*(u)$ satisfies Assumption \ref{a1}. If case (a) holds, then $K^*(u) = 2d\Lambda_1 \indicator(||u||\leq 2L)$. If case (b) holds, then, $K^*(u) = 2d[\Lambda_1 \indicator(||u||\leq 2L) + \left( ||u|| - L\right)^{-\nu} \indicator(||u||> 2L)]$. In both cases $K^*$ is bounded and integrable and therefore satisfies Assumption \ref{a1}. 

Define
\[
a_N := \left(\frac{\ln N}{N h_N^{d_X}}\right)^{1/2}.
\]

\begin{theorem}[\textsc{Weak Uniform Convergence}]
\label{thm:UniformConvergence}
Under Assumptions \ref{a1}, \ref{a2}, \ref{a3}, and the bandwidth conditions $\max\left\{\min\left\{(a_N h_N^{2d_X})^{-\frac{1}{s-1}}, [N^2(\ln (\ln N))^2 \ln N]^{\frac{1}{s}}\right\}, a_N^{-\frac{1}{s-1}}\right\} \ll \min\left\{a_N^{-1}, \frac{N}{\ln N} h_N^{\frac{3}{2}d_X} \right\}$ and $\frac{N }{\ln N}h_N^{d_X} \rightarrow \infty$,
we have for any $q > 0$, $c_N = N^q$,
\[
\sup_{||w|| \leq c_N} \left|\hat{\Psi}_N(w) - \EE\hat{\Psi}_N(w)\right| = O_P(a_N).
\]
\end{theorem}

This theorem establishes uniform convergence of $\hat{\Psi}_N(w)$ to its mean in probability over an expanding set with radius growing at a polynomial rate.

In the proof, we decompose $\hat{\Psi}_N(w)$ into two parts
\[
\hat{\Psi}_N(w) = \tilde{\Psi}_N(w) + R_N(w),
\]
in which
$\tilde{\Psi}_N(w) = \frac{1}{N(N-1)}\sum_{1\leq i \neq j\leq N} Y_{ij}\cdot \indicator\left(|Y_{ij}| < \tau_N\right) K_{ij, N}
$
is a truncated version of $\hat{\Psi}_N(w)$ with a carefully chosen threshold parameter $\tau_N$ and $R_N(w)$ is a residual. The boundedness induced by this truncation is technically convenient as it facilitates the application of various concentration inequalities. To establish concentration of $\tilde{\Psi}_N$, we apply Bernstein inequality to its H\'{a}jek Projection (i.e., to the first-order terms in the Hoeffding decomposition) and apply \cite{arcones1993}'s concentration inequalities for degenerate U-statistics to the second-order terms in the Hoeffding decompositon. Both these bounds requires the truncation threshold to be small enough. To bound the magnitude of the residual $R_N$, we can either apply a triangular inequality to bound the sup of its first moment or use the Borel-Cantelli Lemma to bound its probability of being nonzero. Both these bounds requires the truncation threshold to be large. 

A proper truncation threshold satisfying both requirements exists only if the bandwidth sequence satisfies the condition 
\[
\max\left\{\min\left\{(a_N h_N^{2d_X})^{-\frac{1}{s-1}}, [N^2(\ln (\ln N))^2 \ln N]^{\frac{1}{s}}\right\}, a_N^{-\frac{1}{s-1}}\right\} \ll \min\left\{a_N^{-1}, \frac{N}{\ln N} h_N^{\frac{3}{2}d_X} \right\}. 
\]
The complicated form of this condition is technical in nature. When all (conditional) moments of $Y_{12}$ are bounded, such that $s = \infty$ (of Assumption \ref{a3} above), this condition simplifies to $\frac{N}{\ln N} h_N^{\frac{3}{2} d_X} \gg 1$.

In order to state the weak uniform convergence result for the kernel regression estimator $\hat{g}_N$, we need  additional smoothness assumptions on the kernel. As in other applications of kernel estimation, these assumptions are employed for bias reduction purpose. 
\begin{assumption}[\textsc{Kernel, Part B}]
\label{a:kernel2}
\begin{align*}
\int_{\RR^{d_W}} w_1^{l_1} \cdots w_{d_W}^{l_{d_W}} K(w) dw = 
\left\{
\begin{array}{ll}
	1, & \text{if } l_1 = \cdots = l_{d_W} = 0\\
	0, & \text{if } (l_1, \ldots, l_{d_W})' \in \ZZ_+^{d_W} \text{ and } l_1 + \cdots + l_{d_X} < \beta
\end{array}
\right.
\end{align*} 
\end{assumption}
We can now give a uniform convergence result for the NW regression estimator under dyadic dependence over a sequence of expanding sets.
\begin{theorem}
\label{thm:UniformConvergenceRegression}
Suppose $f_W, g \in \Sigma(\beta, L)$ and
$\delta_N = \inf_{||w||\leq C_N} f_W(w)>0$, $\delta_N^{-1} a_N^*\rightarrow 0$ where $a_N^* := \left(\frac{\ln N}{N h_N^{d_X}}\right)^{1/2} + h_N^{\beta}$. Under the Assumptions of Theorem,
\ref{thm:UniformConvergence}
and Assumption \ref{a:kernel2} 
\begin{align*}
\sup_{||w||\leq C_N} |\hat{g}_N(w) - g(w)| & = O_p(\delta_N^{-1} a_N^*).
\end{align*}
The optimal convergence rate is 
\begin{align*}
\sup_{||w||\leq C_N} |\hat{g}_N(w) - g(w)| & = O_p\left(\delta_N^{-1} \left(\frac{\ln N}{N }\right)^{\frac{\beta}{2\beta + d_X}}\right).
\end{align*}
\end{theorem}
As in the iid case, the KW estimator achieves the optimal rate suggested by Theorem \ref{thm:Minimax} for a compact set with $C_N = C$. If we look at a sequence of expanding sets approaching the entire space $\RR^{d_W}$, then there is an additional penalty term $\delta_N$ due to the presence of the denominator $f_W(w)$.

\bibliography{DyadicNPRegression.bib}

\section*{Appendix}

All notation is as established in the main text unless noted otherwise. Equation numbering continues in sequence with that of the main text.

\subsection*{Proof of Theorem \ref{thm:Minimax}}
Our method of proof follows the general approach outlined in Chapter 2 of \cite{tsybakov2008}. To prove part (i) we use Le Cam's two-point method to find a lower risk bound for estimation of the regression function at a point. To prove statement (ii), which involves the infinity-norm metric, we use Fano's method.

\subsubsection*{Proof of statement (i)}
Our proof of statement (i) essentially involves checking the conditions, as specially formulated for our dyadic regression problem, of Theorem 2.3 of \cite[]{tsybakov2008}. 

For $k = 0,1$, let $P_{kN}$ be a probability measure for the observed data ${\{(X_i', Y_{ij})\}}_{1 \leq i \neq j \leq N}$ with regression function $g_{kN}$. The general reduction scheme outlined in Section 2.2 of \cite{tsybakov2008}, as well as his Theorems 2.1 and 2.2, imply that our Theorem \ref{thm:Minimax} will hold if we can construct two sequences of hypotheses $g_{0N}, g_{1N}$ such that
\begin{enumerate}[label=(\alph*)]
    \item the regression functions $g_{0N}, g_{1N}$ are in the H\"{o}lder class $\Sigma(\beta, L)$;
    
    \item $d\left(\theta_1, \theta_0\right) = |g_{1N}(w) - g_{0N}(w)| \geq  2A\psi_N$ with ${\psi_N} = N^{-\frac{\beta}{2\beta + d_X}}$ and $\theta_0=g_{0N}(w)$ and $\theta_1=g_{1N}(w)$ for some fixed $w \in \mathbb{X} \times \mathbb{X}$;
    
    \item the Kullback-Leibler divergence of $P_{0N}$ from $P_{1N}$ is bounded: $\mathrm{KL}(P_{0N}, P_{1N}) \leq \alpha < \infty$.
\end{enumerate}

The ``trick" of the proof is choosing these two sequences of hypotheses appropriately. Letting $w = (x_{10}, x_{20})$ we choose the sequences:
\begin{align*}
    g_{0N}(x_1, x_2) & \equiv 0 \\
    g_{1N}(x_1, x_2) & = \frac{Lh_N^{\beta}}{2}\left[K\left(\frac{x_1 - x_{10}}{h_N}\right) + K\left(\frac{x_1 - x_{20}}{h_N}\right)+ K\left(\frac{x_2 - x_{10}}{h_N}\right) + K\left(\frac{x_2 - x_{20}}{h_N}\right)\right]
\end{align*}
where $h_N = c_0 N^{-\frac{1}{2\beta + d_X}}$ and the function $K:\RR^{d_W}\rightarrow [0, \infty)$ satisfies
\begin{align}
K \in \Sigma\left(\beta, 1/2\right)\cap C^{\infty}(\RR^{d_X}) \text{ and } K(x) > 0 \Longleftrightarrow ||x||_{\infty} \in (-1/2, 1/2). \label{conditionK}
\end{align}
There exist functions $K$ satisfying this condition. For example, for a sufficiently small $a > 0$, we can take 
\begin{align*}
K(x) = \Pi_{i=1}^{d_X} \lambda(x_i),\ \text{ where }
\lambda(u) = a\eta(2u)\ \text{ and }
\eta(u) = \exp\left(-\frac{1}{1-u^2}\right)\indicator(|u| \leq 1).
\end{align*}
See also Equation (2.34) in \cite{tsybakov2008}.

We verify conditions (a), (b) and (c) in sequence.

\subsubsection*{Verification of (a) $g_{0N}, g_{1N} \in \Sigma(\beta, L)$}
For $s = (\underbrace{s_1, \ldots, s_{d_X}}_{\mathcal{S}_1}, \underbrace{s_{d_{X+1}}, \ldots, s_{2d_X}}_{\mathcal{S}_2})$ with $|s| = \floor{\beta}$, $w = (x_1, x_2)$ and $w' = (x'_1, x'_2)$, the $s^{th}$ order derivative of $g_{1N}$ is
\begin{align*}
D^s g_{1N}(w) & = L h_N^{\beta} \left[D^s K\left(\frac{x_1 - x_{10}}{h_N}\right) + D^s K\left(\frac{x_1 - x_{20}}{h_N}\right) + D^s K\left(\frac{x_2 - x_{10}}{h_N}\right) + D^s K\left(\frac{x_2 - x_{20}}{h_N}\right)\right]\\
& = \left\{
\begin{array}{ll}
	0 & \text{if } |\mathcal{S}_1| \notin \{0, |s|\}\\
	\frac{L h_N^{\beta-\floor{\beta}}}{2} \left[D^{\mathcal{S}_1} K\left(\frac{x_1 - x_{10}}{h_N}\right) + D^{\mathcal{S}_1} K\left(\frac{x_1 - x_{20}}{h_N}\right)\right] & \text{if } |\mathcal{S}_1| = |s|\\
	\frac{L h_N^{\beta-\floor{\beta}}}{2} \left[D^{\mathcal{S}_2} K\left(\frac{x_2 - x_{10}}{h_N}\right) + D^{\mathcal{S}_2} K\left(\frac{x_2 - x_{20}}{h_N}\right)\right] & \text{if } |\mathcal{S}_1| = 0\\	
\end{array}\right..
\end{align*}
Therefore, if $|\mathcal{S}_1| \notin \{0, |s|\}$, then $\left|D^s g_{1N}(w) - D^s g_{1N}(w')\right| = 0$; if $|\mathcal{S}_1| = |s|$, then 
\begin{align*}
& \left|D^s g_{1N}(w) - D^s g_{1N}(w')\right|\\
& = 
\frac{L h_N^{\beta-\floor{\beta}}}{2} \left[\left|D^{\mathcal{S}_1} K\left(\frac{x_1 - x_{10}}{h_N}\right) - D^{\mathcal{S}_1} K\left(\frac{x'_1 - x_{10}}{h_N}\right)\right| + \left|D^{\mathcal{S}_1} K\left(\frac{x_1 - x_{20}}{h_N}\right) - D^{\mathcal{S}_1} K\left(\frac{x'_1 - x_{20}}{h_N}\right)\right|\right]\\
& \leq L ||x_1 - x'_1||_{\infty}^{\beta - \floor{\beta}}\\
& \leq L||w-w'||_{\infty}^{\beta - \floor{\beta}};
\end{align*}
and, finally, if $|\mathcal{S}_1| = 0$, then 
\begin{align*}
& \left|D^s g_{1N}(w) - D^s g_{1N}(w')\right|\\
& = 
\frac{L h_N^{\beta-\floor{\beta}}}{2} \left[\left|D^{\mathcal{S}_2} K\left(\frac{x_2 - x_{10}}{h_N}\right) - D^{\mathcal{S}_2} K\left(\frac{x'_2 - x_{10}}{h_N}\right)\right| + \left|D^{\mathcal{S}_2} K\left(\frac{x_2 - x_{20}}{h_N}\right) - D^{\mathcal{S}_2} K\left(\frac{x'_2 - x_{20}}{h_N}\right)\right|\right]\\
& \leq L ||x_2 - x'_2||_{\infty}^{\beta - \floor{\beta}} \\
& \leq L||w-w'||_{\infty}^{\beta - \floor{\beta}}.
\end{align*}
Hence $g_{1N} \in \Sigma(\beta, L)$. We also have that $g_{0N} \in \Sigma(\beta, L)$ by inspection.

\subsubsection*{Verification of (b): $d\left(\theta\left(P_{0N}\right), \theta\left(P_{1N}\right)\right) = |g_{1N}(w) - g_{0N}(w)| \geq  2A\psi_N$ with $\psi_N = N^{-\frac{\beta}{2\beta + d}}$}
Here we check that our hypotheses are $2s$-separated. We have that
\begin{align*}
|g_{1N}(w) - g_{0N}(w)|
& = \frac{Lh_N^{\beta}}{2} \left[2K\left(0\right) + K\left(\frac{x_{10} - x_{20}}{h_N}\right) + K\left(\frac{x_{20} - x_{10}}{h_N}\right)\right] \geq 2 Lh_N^{\beta} K\left(0\right) \\
& = LK\left(0\right) c_0^{\beta} \psi_N,
\end{align*}
and hence condition (b) holds with $A = \frac{L K \left(0\right) c_0^{\beta}}{2}$.

\subsubsection*{Verification of (c):  $\mathrm{KL}(P_{0N}, P_{1N}) \leq \alpha < \infty$}
This condition allows for the application of part (iii) of Theorem 2.2 in \cite{tsybakov2008}. We begin by establishing some helpful notation. Let $\bY = [Y_{ij}]_{1 \leq i,j \leq N}$ be the $N \times N$ adjacency matrix; $\bG_k=[g_{kN}(W_{ij})]_{1 \leq i,j \leq N}$ for $k=0,1$ the associated matrices of regression functions for the two sequences of hypotheses; and $\bV=[V_{ij}]_{1 \leq i,j \leq N}$ the corresponding matrix of dyadic-specific disturbances. Note the diagonals of each of these matrices consist of ``structural" zeros. Further let $\bU=[U_i]_{1 \leq i \leq N}$ be the $N \times 1$ vector of agent-specific disturbances.
Finally let $\bK$ be the $N \times 1$ vector with $i^{th}$ element $\frac{Lh_N^{\beta}}{2}\left[K\left(\frac{X_i - x_{10}}{h_N}\right) + K\left(\frac{X_i - x_{20}}{h_N}\right)\right]$.

Let $\iota_J$ denote a $J \times 1$ vector of ones, $\b0_{K,J}$ a $K \times J$ matrix of zeros, and $I_J$ the $J \times J$ identity matrix. We also define the following selection matrices:
\begin{align*}
\cT_{1} = 
\begin{pmatrix}
    \iota_{N-1} & 0           & 0      & \cdots & 0      & 0        \\
    \b0         & \iota_{N-2} & 0      & \cdots & 0      & 0        \\
    \vdots      & \vdots      & \ddots & \vdots & \vdots & \vdots   \\
    0           & 0           & 0      & \cdots & 1      & 0        \\
\end{pmatrix}_{\binom{N}{2} \times N},\,\,    
\cT_{2} = 
\begin{pmatrix}
    \b0_{N-1,1}   & I_{N-1}      \\
    \hdashline
    \b0_{N-2,2}   & I_{N-2}      \\
    \vdots        & \vdots       \\
    \hdashline
    0             & 1            \\
\end{pmatrix}_{\binom{N}{2} \times N},
\end{align*}
from which we form $\cT=\cT_1+\cT_2$ and, finally, $\bT=\iota_2 \otimes \cT$.
Next let $\by=(\vech(\bY')',\vech(\bY)')'$ be the $N(N-1) \times 1$ vectorization of the dyadic outcomes. Similarly let $\bg_{k}$ for $k=0,1$ and $\bv$ be the corresponding vectorizations of, respectively, $\bG_k$ and $\bV$.

Using this notation we can write the $N(N-1) \times 1$ vector of composite regression errors $e_{ij} = U_i + U_j + V_{ij}$ as $\mathbf{e} = \bT \bU + \mathbf{v}$ and its variance covariance matrix as
\begin{align*}
    \Omega = \var\left(\mathbf{e}\right) =
    \bI_{N(N-1)\times N(N-1)} + \bT\bT^T.
\end{align*}
Under $P_{0N}$ we have that
\begin{align*}
    \bg_{0} = 0,\ \by = \mathbf{e},\ \by|\bX \sim \Normal\left(0, \Omega\right).
\end{align*}
While under $P_{1N}$ we instead have that
\begin{align*}
    \bg_{1} = \bT\bK,\ \by = \bT\bK + \mathbf{e},\ \by|\bX \sim \Normal\left(\bT\bK, \Omega\right).
\end{align*}

Let $K_{\mathrm{max}} = \max_{u} K(u)$ and recall that $h_N = c_0 N^{-\frac{1}{2\beta + d_X}}$. We can now evaluate the KL divergence as follows:
\begin{align}
\operatorname{KL}\left(P_{0N}, P_{1N}\right) & = \int \log \frac{\mathrm{d}P_{0N}}{\mathrm{d}P_{1N}} \mathrm{d}P_{0N}\\ \nonumber
& = \int \log \frac{p_{0N}(\by|\bX)}{p_{1N}(\by|\bX)} \mathrm{d}P_{0N}\\ \nonumber
& = -\frac{1}{2}\int \by^\top \Omega^{-1} \by - (\by-\bg_{1})^\top\Omega^{-1}(\by-\bg_{1})\mathrm{d}P_{0N}\\ \nonumber
& = \frac{1}{2}\int \bg_{1}^\top\Omega^{-1} \bg_{1} \mathrm{d}P_{0N}\\ \nonumber
& = \frac{1}{2} \EE_{P_{0N}}\left[\bK^\top \bT^\top (\bI+\bT\bT^\top)^{-1}\bT\bK\right]\\ \nonumber
& \leq \frac{1}{2} \EE_{P_{0N}}\left[\bK^\top \bK\right]\\ \nonumber
& \leq \frac{1}{2} L^2 K_{\mathrm{max}}^2 B_3 h_N^{2\beta+d_X} N\\ \nonumber
& = \frac{1}{2} L^2 K_{\mathrm{max}}^2 B_3 c_0^{2\beta+d_X},
\end{align}
for $N$ large enough such that $Nh_{N}^{d_X} \geq 1$ and $LK_{\mathrm{max}}h_{N}^{2\beta}$ bounded above. 

In the derivation above, the third equality follows from the form of the multivariate normal density. The weak inequality in line six holds because
\begin{align*}
\bK^\top \bK - \bK^\top \bT^\top (\bI+\bT\bT^\top)^{-1}\bT\bK & = 
\bK^\top \left[\bI_N - \bT^\top (\bI+\bT\bT^\top)^{-1}\bT\right]\bK\\
& = \bK^\top \left[\bI_N + \bT^\top \bT\right]^{-1}\bK\\
& \geq 0.
\end{align*}
Finally, the weak inequality in line seven holds because, using condition \eqref{conditionK} above,
\begin{align*}
& \EE\left[\left(K\left(\frac{X_i - x_{10}}{h_N}\right) + K\left(\frac{X_i - x_{20}}{h_N}\right)\right)^2\right]\\
& \leq 2\EE\left[\left(K\left(\frac{X_i - x_{10}}{h_N}\right)\right)^2 + \left(K\left(\frac{X_i - x_{20}}{h_N}\right)\right)^2\right]\\
& = 2\int \left(K\left(\frac{x - x_{10}}{h_N}\right)\right)^2 + \left(K\left(\frac{x - x_{20}}{h_N}\right)\right)^2 \mathrm{d}F(x)\\
& \leq 2 K_{\mathrm{max}}^2 \int \indicator\left(\left|\frac{x - x_{10}}{h_N}\right| \leq \frac{1}{2}\right) + \indicator\left(\left|\frac{x - x_{20}}{h_N}\right| \leq \frac{1}{2}\right) \mathrm{d}F(x)\\
& = 2 K_{\mathrm{max}}^2 h_N^{d_X}\left[\int \indicator\left(\left|u\right| \leq \frac{1}{2}\right) [f(x_{10} + h_N u) + f(x_{20} + h_N u)] \mathrm{d}u\right] \\
& \leq 4 h_N^{d_X} B_3 K_{\mathrm{max}}^2,
\end{align*}
and where it is also helpful to remind oneself of the definition of $\bK$ given earlier.

If we take $c_0 = \left(\frac{2\alpha}{L^2 K_{max}^2 B_3}\right)^{\frac{1}{2\beta+d_X}}$, then we obtain $\operatorname{KL}\left(P_{0N}, P_{1N}\right) \leq \alpha$. This result, and condition (b) above, gives -- invoking
equations (2.7) and (2.9) on p. 29 of \cite{tsybakov2008} as well as part (iii) of his Theorem 2.2:
\begin{equation*}
    \underset{\hat{g}_{N}}{\inf}\underset{g\in\Sigma\left(\beta,L\right)}{\sup}\mathbb{E}_{g}\left[1\left(|g_{1N}(w)-g_{0N}(w)|\geq A\psi_{N}\right)\right]\geq\max\left(\frac{1}{4}\exp\left(-\alpha\right),\frac{1-\sqrt{\frac{\alpha}{2}}}{2}\right)
\end{equation*}
for $N$ large enough. Some rearrangement and the Markov Inequality then yield
\begin{equation*}
    \underset{\hat{g}_{N}}{\inf}\underset{g\in\Sigma\left(\beta,L\right)}{\sup}\mathbb{E}_{g}\left[N^{\frac{2\beta}{2\beta+d_{X}}}\left(g_{1N}(w)-g_{0N}(w)\right)^{2}\right]\geq A^{2}\max\left(\frac{1}{4}\exp\left(-\alpha\right),\frac{1-\sqrt{\frac{\alpha}{2}}}{2}\right).
\end{equation*}
Since the constant to the right of the inequality only depends on $\beta$ and $L$ part (i) of the Theorem follows after taking the limit inferior of the expression above as $N \rightarrow \infty$.

\subsubsection*{Proof of statement (ii)}
Again let $P_{kN}$ be the probability measure of the observed data $(X_i, Y_{ij}, 1 \leq i \neq j \leq N)$ with the regression function $g_{kN}$. Theorem 2.5 of \cite{tsybakov2008} implies that part (ii) will hold if we can construct sequences of hypotheses $P_{0N}, P_{1N}, \ldots, P_{M_N N}$ such that
\begin{enumerate}[label=(\alph*)]
\item $g_{0N}, g_{kN} \in \Sigma(\beta, L)$, $k = 1, \ldots, M_N$;
\item $d\left(\theta_k, \theta_l\right) = ||g_{kN} - g_{lN}||_{\infty} \geq  2A\psi_N$, ${\psi_N} = \left(\frac{N}{\ln N}\right)^{-\frac{\beta}{2\beta + d}}$ and $\theta_k=g_{kN}$ and $\theta_l=g_{lN}$ for $k \neq l$ and $k,l=1,\dots,M_N$;
\item $\frac{1}{M_N}\sum_{k = 1}^{M_N}KL(P_{kN}, P_{0N}) \leq \alpha \ln M_N$.
\end{enumerate}
Define the hypotheses:
\begin{align*}
g_{0N}: & (x_1, x_2) \rightarrow 0\\
g_{kN}: & (x_1, x_2) \rightarrow Lh_N^{\beta}\left[K\left(\frac{x_1 - x_{kN}}{h_N}\right) + K\left(\frac{x_2 - x_{kN}}{h_N}\right)\right]
\end{align*}
where $k \in \cI_N = \{1, 2, \ldots, m_N\}^{d_X}$, $h_N = c_0 \left(\frac{N}{\ln N}\right)^{-\frac{1}{2\beta + d_X}}$, $m_N = \ceil{h_N^{-1}}$, $M_N = |\cI_N| = m_N^{d_X}$, and for $k = (k_1, k_2, \ldots, k_d)$, $x_{kN} = \left(\frac{k_1 - 1/2}{m_N}, \frac{k_2 - 1/2}{m_N}, \ldots, \frac{k_d - 1/2}{m_N}\right)$,
 the function $K:\RR^{d_X}\rightarrow [0, \infty)$ satisfies (\ref{conditionK}). Notice the supports of these functions for the same $N$ are disjoint.
The results follows by verifying conditions (a), (b) and (c). We have already shown that condition (a) holds in the proof of part (i). The condition (b) holds with $A = LK(0) c_0^{\beta}$ because
\begin{align*}
||g_{kN} - g_{lN}||_{\infty}
\geq |g_{kN}(x_{kN}, x_{kN}) - g_{lN}(x_{kN}, x_{kN})|
= 2Lh_N^{\beta}K(0) = 2LK(0) c_0^{\beta} \psi_N.
\end{align*}
To verify condition (c) we evaluate the KL-divergence:
\begin{align*}
\frac{1}{M_N}\sum_{k\in \cI_N} \mathrm{KL}(P_{kN}, P_{0N})
& \leq 
\frac{1}{M_N}\sum_{k\in \cI_N} \frac{1}{2} \EE_{P_{0N}}\left[\bK_k^\top \bK_k\right]\\
& \leq 
\frac{1}{M_N}\sum_{k\in \cI_N} 2 L^2 h_N^{2\beta} 
K_{\mathrm{max}}^2 \sum_{i=1}^N \int \indicator\left(\left|\frac{x_i - x_{kN}}{h_N}\right| \leq \frac{1}{2}\right) \mathrm{d}F(x_i)\\
& = 
\frac{1}{M_N} 2 L^2 h_N^{2\beta} 
K_{\mathrm{max}}^2 \sum_{i=1}^N \int \sum_{k\in \cI_N} \indicator\left(\left|\frac{x_i - x_{kN}}{h_N}\right| \leq \frac{1}{2}\right) \mathrm{d}F(x_i)\\
& \leq
2 L^2 h_N^{2\beta + d_X} K_{\mathrm{max}}^2 N\\
& = 2 L^2 K_{\mathrm{max}}^2 c_0^{2\beta + d_X}\ln N.
\end{align*}
The first and second line are proved in part (i). The fourth line use the fact that the support of functions $g_{kN}, k\in \cI_N$ are disjoint and $\sum_{k\in \cI_N} \indicator\left(\left|\frac{x_i - x_{kN}}{h_N}\right| \leq \frac{1}{2}\right) \leq 1$.
We have
$\ln M_N = \ln (m_N^{d_X}) \geq \frac{d_X}{2\beta + d_X} \ln\left(\frac{N}{\ln N}\right) - d_X \ln c_0 \geq \frac{d_X}{2\beta + d_X +1} \ln N$ for sufficiently large $N$. The condition is thus satisfied with sufficiently large $c_0$. The result follows from Theorem 2.5 of \cite{tsybakov2008}.

\subsection*{Proof of Theorem \ref{thm:VarianceBound}}
Applying the variance operator to $\hat{\Psi}(w)$ yields
\begin{equation*}
\mathbb{V}\left(\hat{\Psi}(w)\right) 
=  \frac{4}{N}\frac{N-2}{N-1} V_{N,1} + \binom{N}{2}^{-1} V_{N,2} 
\label{eq:variance}
\end{equation*}
where, starting with the second term, 
\begin{align*}
V_{N,2} & = \mathbb{V}\left(\frac{1}{2}\left[Y_{12} K_{12} + Y_{21} K_{21}\right]\right) \leq \mathbb{V}\left(Y_{12} K_{12}\right) \leq \EE\left(Y_{12}^2 K^2_{12}\right)\\
& = h_N^{-4d_X}\int \EE\left[Y_{12}^2 | (X_1, X_2) = (x_1, x_2)\right] K^2\left(\frac{x-x_1}{h_N}, \frac{x-x_2}{h_N}\right) f(x_1) f(x_2) \mathrm{d}x_1\mathrm{d}x_2\\
& = h_N^{-2d_X} \int \EE\left[Y_{12}^2 | (X_1, X_2) = (x - h_Ns_1, x' - h_Ns_2)\right]  f(x - h_N s_1) f(x' - h_N s_2) K^2\left(s_1, s_2\right)\mathrm{d}s_1\mathrm{d}s_2\\
& \leq h_N^{-2d_X} B_4 K_{\mathrm{max}} B_1.
\end{align*}
Next, consider the first term. We get that
\begin{align*}
V_{N,1} & = \mathbb{C}\left(
\frac{1}{2}
\left(Y_{12} K_{12} + Y_{21} K_{21}\right),
\frac{1}{2}
\left(Y_{13} K_{13} + Y_{31} K_{31}\right)
\right)\\
& = \mathbb{V}\left(
\EE\left[
\frac{1}{2}
\left(Y_{12} K_{12} + Y_{21} K_{21}\right) \Big|X_1, U_1\right]
\right)\\
& \leq \frac{1}{2}\var\left(
\EE\left(
Y_{12} K_{12} \Big|X_1, U_1\right)
\right)
+
\frac{1}{2}
\var\left(
\EE\left(
Y_{21} K_{21} \Big|X_1, U_1\right)
\right)
\\
& \leq 
\frac{1}{2}\EE\left(Y_{12} K_{12} Y_{13} K_{13}\right)
+ 
\frac{1}{2}\EE\left(Y_{21} K_{21} Y_{31} K_{31}\right)\\
& = \frac{1}{2} h_N^{-4d_X}\int \EE\left(Y_{12}Y_{13}|(X_1, X_2, X_3) = (x_1, x_2, x_3)\right) \\
& \qquad \qquad \qquad \cdot K\left(\frac{x - x_1}{h_N}, \frac{x' - x_2}{h_N}\right) K\left(\frac{x - x_1}{h_N}, \frac{x' - x_3}{h_N}\right) f(x_1) f(x_2) f(x_3)\mathrm{d}x_1 \mathrm{d}x_2 \mathrm{d}x_3\\
& \quad + \frac{1}{2} h_N^{-4d_X}\int \EE\left(Y_{21}Y_{31}|(X_1, X_2, X_3) = (x_1, x_2, x_3)\right) \\
& \qquad \qquad \qquad \cdot K\left(\frac{x - x_2}{h_N}, \frac{x' - x_1}{h_N}\right) K\left(\frac{x - x_3}{h_N}, \frac{x' - x_1}{h_N}\right) f(x_1) f(x_2) f(x_3)\mathrm{d}x_1 \mathrm{d}x_2 \mathrm{d}x_3\\
& = h_N^{-d_X} \frac{1}{2}\int \EE\left(Y_{12}Y_{13}|(X_1, X_2, X_3) = (x - h_N s_1, x' - h_N s_2, x' - h_N s_3)\right) \\
& \qquad \qquad \cdot f(x - h_N s_1) f(x'- h_N s_2) f(x'- h_N s_3)K\left(s_1, s_2\right) K\left(s_1, s_3\right) \mathrm{d}s_1 \mathrm{d}s_2 \mathrm{d}s_3\\
& \quad + h_N^{-d_X} \frac{1}{2}\int \EE\left(Y_{21}Y_{31}|(X_1, X_2, X_3) = (x' - h_Ns_1, x - h_Ns_2, x - h_Ns_3) \right) \\
& \qquad \qquad \cdot f(x' - h_N s_1) f(x- h_N s_2) f(x- h_N s_3)K\left(s_1, s_2\right) K\left(s_1, s_3\right) \mathrm{d}s_1 \mathrm{d}s_2 \mathrm{d}s_3\\
& \leq h_N^{-d_X} B_5 \int |K\left(s_1, s_2\right)| |K\left(s_1, s_3\right)| \mathrm{d}s_1 \mathrm{d}s_2 \mathrm{d}s_3\\
& \leq h_N^{-d_X} B_5 B_2 B_1. \numberthis \label{vn1UpperBound}
\end{align*}
These two bounds imply the variance bound
\begin{align*}
\mathbb{V}\left(\hat{\Psi}(w)\right) 
& \leq \binom{N}{2}^{-1} h_N^{-2d_X} B_4 K_{\mathrm{max}} B_1 + \frac{4(N-2)}{N(N-1)} h_N^{-d_X} B_5 B_2 B_1\\
& = N^{-1}h_N^{-d_X} \left[\frac{N-2}{N-1}4B_5 B_2 B_1 + N^{-1}h_N^{-d_X}\frac{4N}{N-1}B_4 K_{\mathrm{max}} B_1\right],
\end{align*}
which, in turn, implies that for $M_0 = 4B_5 B_2 B_1 + 1$ and sufficiently large $N$,
$\mathbb{V}\left(\hat{\Psi}(w)\right) \leq \frac{M_0}{N h_N^{d_X}}$ for all $w\in \RR^{d_W}$ as claimed.

\subsection*{Proof of Theorem \ref{thm:UniformConvergence}}
For $\tau_N$ a sequence of positive truncation parameters we consider the sum
\begin{align*}
\tilde{\Psi}_N(w) = \frac{1}{\binom{N}{2}}\sum_{1\leq i<j\leq N} \frac{1}{2} 
                    & \left[ Y_{ij}\cdot \indicator\left(|Y_{ij}| < \tau_N\right) \frac{1}{h_N^{d_W}} K\left(\frac{w-W_{ij}}{h_N}\right) \right. \\
                    & \left. + Y_{ji}\cdot \indicator\left(|Y_{ji}| < \tau_N\right) \frac{1}{h_N^{d_W}} K\left(\frac{w-W_{ji}}{h_N}\right)\right].
\end{align*}
We will use $\tilde{Z}_{N, ij}$ to denote the summands in the above expression in what follows. The Hoeffding decomposition of this $U$-like statistic is
\begin{align*}
\tilde{\Psi}(w) = \EE \tilde{\Psi}(w) + \underbrace{\frac{2}{N}\sum_{i=1}^N \bar{Z}_{N,i}}_{T_{N, 1}(w)} + \underbrace{\frac{1}{\binom{N}{2}} \sum_{1 \leq i < j \leq N} \breve{Z}_{N, ij}}_{T_{N,2}(w)},
\end{align*}
where
\begin{align*}
\bar{Z}_{N,i} & = \EE\left[\tilde{Z}_{N, ij} \Big |X_i, U_i\right] - \EE \tilde{Z}_{N, ij}\\
\breve{Z}_{N, ij} & = \tilde{Z}_{N, ij} - \EE\left[\tilde{Z}_{N, ij} \Big |X_i, U_i\right] - \EE\left[\tilde{Z}_{N, ij} \Big |X_j, U_j\right] + \EE \tilde{Z}_{N, ij}.
\end{align*}
Notice that $T_{N, 1}(w)$ is an average of $N$ iid mean-zero random variables while $T_{N, 2}(w)$ is a degenerate second-order $U$-like statistic. 

To proceed further we require the following Lemma.
\begin{lemma}
\label{lemma:concentration}
Under Assumptions \ref{a1} and \ref{a2}, for any $\alpha > 0$, there exists constant $M_{\alpha}$ such that
\begin{enumerate}[label=(\roman*)]
\item if $\tau_N \ll a_N^{-1}$, then $\sup_{w\in \RR^{d_W}}P\left(|T_{N, 1}(w)| > M_{\alpha} a_N\right) = O\left(N^{-\alpha}\right)$;
\item if $\tau_N \ll N h^{\frac{3}{2}d_X} /\ln N$ and $a_N = o(1)$, then $\sup_{w\in \RR^{d_W}}P\left(|T_{N, 2}(w)| > M_{\alpha} a_N\right) = O\left(N^{-\alpha}\right)$;
\item if for some $s>1$, $\sup_{x_1, x_2 \in \RR^{d_X}} \EE\left[|Y_{12}|^s\big|(X_1, X_2)=(x_1, x_2)\right] \cdot f(x_1, x_2) \leq B_{4,s} < \infty$ and $\tau_N \gg a_N ^{-\frac{1}{s-1}}$, then $\sup_{w\in \RR^{d_W}}\left|\EE\left(\hat{\Phi}(w) - \tilde{\Phi}(w)\right)\right| = o\left(a_N\right)$;
\item if for some $s>1$, $\EE\left|Y_{12}\right|^s \leq B_{6, s}$ and $\tau_N \gg (a_N h_N^{2d_X})^{-\frac{1}{s-1}}$, then $\sup_{w\in \RR^{d_W}}\left|\hat{\Phi}_N(w) - \tilde{\Phi}_N(w)\right| = o_P\left(a_N\right)$;
\item if for some $s>2$, $\tau_N = \left(N^2 \phi_N\right)^{\frac{1}{s}}$ where $\phi_N = (\ln (\ln N))^2 \ln N$, and $\EE\left|Y_{12}\right|^s \leq B_{6, s}$, then $P(\hat{\Phi}_N = \tilde{\Phi}_N) = P\left(\hat{\Phi}_N(w) = \tilde{\Phi}_N(w), \forall w \in \RR^{2d_X}\right) \rightarrow 1$ as $N \rightarrow \infty$.
\end{enumerate}
\end{lemma}

The proof of the above Lemma may be found below. The bandwidth conditions stated in the hypotheses of Theorem \ref{thm:UniformConvergence} ensure that we can pick truncation thresholds $\tau_N$ which satisfy the following conditions
\begin{enumerate}
\item $\tau_N \ll a_N^{-1}$;
\item $\tau_N \ll \frac{N}{\ln N} h_N^{\frac{3}{2}d_X} $;
\item $\tau_N \gg a_N^{-\frac{1}{s-1}}$;
\item $\tau_N \gg \left(N^2 \phi_N\right)^{\frac{1}{s}}$ or $\tau_N \gg (a_N h_N^{2d_X})^{-\frac{1}{s-1}}$.
\end{enumerate} 
These conditions allow for the application of Lemma \ref{lemma:concentration}.
Denote $R_N(w) := \hat{\Psi}_N(w) - \tilde{\Psi}_N(w)$. For any set $\cC_N \subset \RR^{2d}$,
\begin{align*}
& P\left(\sup_{w \in \cC_N} \left|\hat{\Psi}_N(w) - \EE\hat{\Psi}_N(w)\right| > 8 M a_N\right)\\
& = P\left(\sup_{w \in \cC_N} \left|\tilde{\Psi}_N(w) - \EE\tilde{\Psi}_N(w) + R_N(w) - \EE R_N(w)\right| > 8 M a_N\right)\\
& \leq P\left(\sup_{w \in \cC_N} \left|\tilde{\Psi}_N(w) - \EE\tilde{\Psi}_N(w)\right| > 6 M a_N\right) + P\left(\sup_{w \in \cC_N} \left|R_N(w) - \EE R_N(w)\right| > 2 M a_N\right). \numberthis \label{ineq:TwoPiecesByTruncation}
\end{align*}
The second term in inequality (\ref{ineq:TwoPiecesByTruncation}) converges to zero because
\begin{align*}
& P\left(\sup_{w \in \cC_N} \left|R_N(w) - \EE R_N(w)\right| > 2 M a_N\right) \\& \leq P\left(\sup_{w \in \RR^{d_W}} \left|R_N(w) - \EE R_N(w)\right| > 2M a_N\right)\\
& \leq P\left(\sup_{w \in \RR^{d_W}} \left|R_N(w)\right| > M a_N\right) + \indicator\left(\sup_{w \in \RR^{d_W}} \left|\EE R_N(w)\right| > M a_N\right)\numberthis \label{ineq:LargeValuePart}\\
& = o(1). 
\end{align*}
The last line holds because 
\begin{align}
& \indicator\left(\sup_{w \in \RR^{d_W}} \left|\EE R_N(w)\right| > M a_N\right) = 0 \qquad \text{for large }N \label{ResidualExpectation}\\
& P\left(\sup_{w \in \RR^{d_W}} \left|R_N(w)\right| > M a_N\right) = o_P(1).\label{Residual}
\end{align}

To see (\ref{ResidualExpectation}), notice part (iii) of Lemma \ref{lemma:concentration} implies that $\sup_{w \in \RR^{d_W}} \left|\EE R_N(w)\right| = o(a_N)$. Hence $\indicator\left(\sup_{w \in \RR^{2d}} \left|\EE R_N(w)\right| > M a_N\right) = 0$ for large $N$. To see (\ref{Residual}), notice the inequality
\begin{align*}
P\left(\sup_{w \in \RR^{d_W}} \left|R_N(w)\right| > M a_N\right) \leq \min\left\{1 - P(\hat{\Phi}_N = \tilde{\Phi}_N), \frac{\EE\sup_{w \in \RR^{d_W}} \left|R_N(w)\right|}{Ma_N}\right\},
\end{align*}
suggests we can bound either term on the right-hand side to bound the term on the left-hand side. The threshold we pick meets the conditions of both parts (iv) and (v) of Lemma \ref{lemma:concentration}, which ensures either $1 - P(\hat{\Phi}_N = \tilde{\Phi}_N) = o(1)$ or $\frac{\EE\sup_{w \in \RR^{d_W}} \left|R_N(w)\right|}{Ma_N} = o(1)$. This implies (\ref{Residual}).

To show the first term in inequality (\ref{ineq:TwoPiecesByTruncation}) converges to zero, we will use a covering argument to reduce finding the supremum over an infinite number points to finding the maximum over a finite number of points. We then invoke point-wise concentration bounds. This part closely follows the argument in \cite{hansen2008}. Cover any compact region $\cC_N \subset \RR^{d_W}$ by finite number of balls of radius $a_N h_N$ centered at grid points in the set $L_N = \{w_{N, 1}, w_{N, 2}, \ldots, w_{N, L_N}\}$ (Here we abuse the notation a bit: $L_N$ is used to refer to both the set and its cardinality). 
Denote the ball $A_{N, j} = \{w \in \RR^{d_W}: ||w - w_{N, j}|| \leq a_N h_N\}$. 
For $N$ large enough such that $a_N < L$ ($L$ is the constant appearing in Assumption \ref{a3}), for any point $w \in A_{N, j}$ within the ball, assumption \ref{a3} (ii) implies
\begin{align}
\left|K\left(\frac{w-W_{ij}}{h}\right) - K\left(\frac{w_{N, j}-W_{ij}}{h}\right)\right| \leq a_N K^*\left(\frac{w_{N, j}-W_{ij}}{h}\right). \label{KernelBound}
\end{align}
Define
\begin{align*}
\breve{\Phi}_N(w) := \frac{1}{N(N-1)}\sum_{1\leq i \neq j\leq N} Y_{ij}\cdot \indicator\left(|Y_{ij}| < \tau_N\right) \frac{1}{h^{d_W}} K^*\left(\frac{w-W_{ij}}{h}\right),
\end{align*}
which is a version of $\tilde{\Phi}(w)$ with $K$ replaced by $K^*$. The bound (\ref{KernelBound}) implies
\begin{align*}
\left|\tilde{\Psi}_N(w) - \tilde{\Psi}_N(w_{N, j})\right| \leq a_N \breve{\Phi}_N(w_{N, j}),
\end{align*}
with 
$|\EE \breve{\Phi}_N(w_{N, j})| \leq B_4^{1/2} B_3^{1/2} \int |K^*(w)|\mathrm{d}w < \infty$.
Next bound the sup within the ball by a value at the center and the sup discrepancy
\begin{align*}
& \sup_{w \in A_{N, j}} \left|\tilde{\Psi}_N(w) - \EE\tilde{\Psi}_N(w)\right| \\
& \leq \left|\tilde{\Psi}_N(w_{N, j}) - \EE\tilde{\Psi}_N(w_{N, j})\right| + \sup_{w \in A_{N, j}} \left|\tilde{\Psi}_N(w) - \tilde{\Psi}_N(w_{N, j})\right| + \sup_{w \in A_{N, j}} \left|\EE\left(\tilde{\Psi}_N(w) - \tilde{\Psi}_N(w_{N, j})\right)\right|\\
& \leq \left|\tilde{\Psi}_N(w_{N, j}) - \EE\tilde{\Psi}_N(w_{N, j})\right| + 
a_N\left[\breve{\Phi}_N(w_{N, j}) + \EE\breve{\Phi}_N(w_{N, j})\right]\\
& \leq \left|\tilde{\Psi}_N(w_{N, j}) - \EE\tilde{\Psi}_N(w_{N, j})\right| + 
a_N\left|\breve{\Phi}_N(w_{N, j}) - \EE\breve{\Phi}_N(w_{N, j})\right| + 2a_N\EE\breve{\Phi}_N(w_{N, j})\\
& \leq \left|\tilde{\Psi}_N(w_{N, j}) - \EE\tilde{\Psi}_N(w_{N, j})\right| + 
\left|\breve{\Phi}_N(w_{N, j}) - \EE\breve{\Phi}_N(w_{N, j})\right| + 2a_N\EE\breve{\Phi}_N(w_{N, j}).
\end{align*}
The last inequality follow because $a_N \leq 1$ for $N$ large enough. For any constant $M \geq B_4^{1/2} B_3^{1/2} \int |K^*(w)|\mathrm{d}w \geq \EE \breve{\Phi}_N(w_{N, j})$,
\begin{align*}
& P\left(\sup_{w \in A_{N, j}} \left|\tilde{\Psi}_N(w) - \EE\tilde{\Psi}_N(w)\right| > 6 M a_N\right) \\
& \leq P\left(\left|\tilde{\Psi}_N(w_{N, j}) - \EE\tilde{\Psi}_N(w_{N, j})\right| + \left|\breve{\Phi}_N(w) - \EE\breve{\Phi}_N(w)\right| + 2a_N\EE\breve{\Phi}_N(w) > 6 M a_N\right)\\
& \leq P\left(\left|\tilde{\Psi}_N(w_{N, j}) - \EE\tilde{\Psi}_N(w_{N, j})\right| > 2 M a_N\right) + P\left(\left|\breve{\Phi}_N(w) - \EE\breve{\Phi}_N(w)\right| > 2 M a_N\right),
\end{align*}
as well as 
\begin{align*}
& P\left(\sup_{w \in \cC_N} \left|\tilde{\Psi}_N(w) - \EE\tilde{\Psi}_N(w)\right| > 6 M a_N\right) \\
& \leq \sum_{j = 1}^{L_N} P\left(\sup_{w \in A_{N, j}} \left|\tilde{\Psi}_N(w) - \EE\tilde{\Psi}_N(w)\right| > 6 M a_N\right)\\
& \leq L_N \max_{j \in \{1, 2, \ldots, L_N\}} P\left(\sup_{w \in A_{N, j}} \left|\tilde{\Psi}_N(w) - \EE\tilde{\Psi}_N(w)\right| > 6 M a_N\right)\\
& \leq L_N \max_{j \in \{1, 2, \ldots, L_N\}} P\left(\left|\tilde{\Psi}_N(w_{N, j}) - \EE\tilde{\Psi}_N(w_{N, j})\right| > 2 M a_N\right)\\
& \qquad + L_N \max_{j \in \{1, 2, \ldots, L_N\}} P\left(\left|\breve{\Phi}_N(w) - \EE\breve{\Phi}_N(w)\right| > 2 M a_N\right). \numberthis \label{6}
\end{align*}
We now bound the two terms in (\ref{6}) using the same argument, as both $K$ and $K^*$ satisfy Assumption \ref{a1}, and this is the only property of the function $K$ or $K^*$ we will use. For any $\alpha > 0$ and $M_{\alpha}$ as in Lemma \ref{lemma:concentration}, for any $w \in \RR^{d_W}$
\begin{align*}
\sup_{w\in \RR^{d_W}} P\left(\left|\tilde{\Psi}_N(w) - \EE\tilde{\Psi}_N(w)\right| > 2 M_{\alpha} a_N\right) 
& = \sup_{w\in \RR^{d_W}}P\left(\left|T_{N, 1}(w) + T_{N, 2}(w)\right| > 2 M_{\alpha} a_N\right)\\
& \leq \sup_{w\in \RR^{d_W}}P\left(\left|T_{N, 1}(w)\right| > M_{\alpha} a_N\right) \\
& \qquad + \sup_{w\in \RR^{d_W}}P\left(\left|T_{N, 2}(w)\right| > M_{\alpha} a_N\right)\\
& = O\left(N^{-\alpha}\right).
\end{align*}
Hence
\begin{align*}
P\left(\sup_{w \in \cC_N} \left|\tilde{\Psi}_N(w) - \EE\tilde{\Psi}_N(w)\right| > 6 M a_N\right) \leq O\left(L_N N^{-\alpha}\right).
\end{align*}
If we take $\cC_N = \{w\in \RR^{d_W}: ||w|| < c_N\}$ where $c_N = N^{q}$, then $\cC_N$ can be covered by $L_N = 2 \left(\frac{c_N}{a_N h_N}\right)^{d_W}$ number of balls with radius $a_N h_N$. Hence we can take $\alpha$ large enough, e.g. $\alpha = (q+\frac{1}{2})d_W+3$, so that 
$O\left(L_N N^{-\alpha}\right) 
= O\left(\left(\frac{c_N}{a_N h_N}\right)^{d_W} N^{-\alpha}\right) 
= O\left(N^{(q+\frac{1}{2})d_W+2-\alpha}\right) 
= O\left(N^{-1}\right) 
= o(1)$. We have therefore shown that 
\begin{align*}
P\left(\sup_{w \in \cC_N} \left|\tilde{\Psi}_N(w) - \EE\tilde{\Psi}_N(w)\right| > 6 M a_N\right) = o(1). \numberthis \label{ineq:BoundedPart}
\end{align*}
Together the two bounds (\ref{ineq:LargeValuePart}) and (\ref{ineq:BoundedPart}) imply that the right-hand side of equality (\ref{ineq:TwoPiecesByTruncation}) is $o(1)$. This is saying for sufficiently large $M < \infty$, we have $P\left(\sup_{w \in \cC_N} \left|\hat{\Psi}_N(w) - \EE\hat{\Psi}_N(w)\right| > 8M a_N\right)=o(1)$, which is sufficient for
\begin{align*}
\sup_{||w|| \leq c_N} \left|\hat{\Psi}_N(w) - \EE\hat{\Psi}_N(w)\right| = O(a_N).
\end{align*}
as required.

\subsection*{Proof of Lemma \ref{lemma:concentration}}

\subsubsection*{Proof of claim (i)}
To prove the first claim of the Lemma we will apply the classic Bernstein's inequality (see equation 2.10 on p. 36 of the textbook \cite{boucheron2013}).
\begin{quote}
Let $Q_1, \ldots, Q_N$ be independent random variables with finite variance such that $Q_i \leq b$ for some $b > 0$ almost surely for all $i<N$. Let $S = \sum_{i=1}^N (Q_i - \EE Q_i)$ and $v = \sum_{i=1}^N \EE\left[Q_i^2\right]$. Then for any $t>0$, 
\[
P(S \geq t) \leq \exp\left(- \frac{t^2}{2(v + bt/3)}\right).
\]
\end{quote}
In order to invoke the inequality, we first show that $Q_i(w) := \tau_N^{-1} h_N^{d_X} \bar{Z}_{N,i}(w)$ is bounded. In the following we will use the abbreviated notation $Q_i$ for $Q_i(w)$. 
Remember $\bar{Z}_{N,i}$ is the mean-normalized version of $\EE\left[\tilde{Z}_{N, ij} \Big |X_i, U_i\right]$. Since
\begin{align*}
\tau_N^{-1} h^{d_X} \left|\EE\left[\tilde{Z}_{N, ij} \Big |X_i, U_i\right]\right|
= & \tau_N^{-1} h^{d_X} \left| \EE \left[ \frac{1}{2}\left[Y_{ij}\cdot \indicator\left(|Y_{ij}| < \tau_N\right) K_{ij} \right. \right. \right. \\
& \left. \left. \left. + Y_{ji}\cdot \indicator\left(|Y_{ji}| < \tau_N\right) K_{ji}\right]\Big |X_i, U_i \right] \right| \\
\leq & h_N^{d_X} \frac{1}{2}\EE\left[|K_{ij}| + |K_{ji}|\Big |X_i, U_i\right]\\
= & h_N^{-d_X} \frac{1}{2}\EE\left[\left|K\left(\frac{w - W_{ij}}{h_N}\right)\right| + \left|K\left(\frac{w - W_{ji}}{h_N}\right)\right| \Big |X_i\right]\\
= & h_N^{-d_X} \frac{1}{2}\int \left[\left|K\left(\frac{x-x_i}{h_N}, \frac{x'-x_j}{h_N}\right)\right| + \left|K\left(\frac{x-x_j}{h_N}, \frac{x'-x_i}{h_N}\right)\right|\right] f(x_j)\mathrm{d}x_j\\
= & \frac{1}{2} \int \left|K\left(\frac{x-x_i}{h_N}, s\right)\right| f(x' - h_N s) + \left|K\left(s, \frac{x'-x_i}{h_N}\right)\right| f(x - h_N s) \mathrm{d}s\\
\leq & B_2 B_3,
\end{align*}
we have $|Q_i| = |\tau_N^{-1} h_N^{d_X} \bar{Z}_{N,i}| < 2 B_2 B_3$. 
Write $P\left(T_{N,1}(w) > M a_N\right)$ in the form suitable for applying Bernstein's inequality
\begin{align*}
P\left(T_{N,1}(w) > M a_N\right) 
& = P\left(\frac{2}{N}\sum_{i=1}^N \bar{Z}_{N,i} > M a_N\right)\\
& = P\left(\sum_{i=1}^N \tau_N^{-1} h_N^{d_X}\bar{Z}_{N,i} > \frac{M}{2}  Nh_N^{d_X}a_N \tau_N^{-1} \right)\\
& = P(S \geq t),
\end{align*}
in which $S = \sum_{i=1}^N Q_i$ and $t = \frac{M}{2}  Nh_N^{d_X}a_N \tau_N^{-1}$. 
Applying Bernstein's inequality gives us $P(S \geq t) \leq \exp\left(- \frac{t^2}{2(v + bt/3)}\right)$ where $v := \sum_{i=1}^N \EE\left[Q_i^2\right]$ and $b = 2B_2B_3$. Since the function $\exp\left(- \frac{t^2}{2(v + bt/3)}\right)$ is increasing in $v$, we have for any $v' > v$
\begin{equation}
P(S \geq t) \leq \exp\left(- \frac{t^2}{2(v' + bt/3)}\right). \label{bernstein}
\end{equation}
The upper bound $v'$ we are going to use is the following one
\begin{align*}
v & = \sum_{i=1}^N \EE\left[Q_i^2\right] = \sum_{i=1}^N \EE\left[\left(\tau_N^{-1} h_N^{d_X} \bar{Z}_{N,i}\right)^2\right] = \tau_N^{-2} h_N^{2d_X} N V_{N, 1} \leq \tau_N^{-2} N h_N^{d_X} B_5 B_2 B_1 := v',
\end{align*}
in which the inequality is an implication of \eqref{vn1UpperBound}. Plugging the expression of $v'$, $t$, $b$, and $a_N$ into the RHS of \eqref{bernstein} gives us
\begin{align*}
\exp\left(- \frac{t^2}{2(v' + bt/3)}\right)
& = 
\exp\left(- \frac{M^2}{8B_5 B_2 B_1 + 8B_2 B_3 M  a_N \tau_N/3} \ln N\right).
\end{align*}
By assumption $a_N\tau_N \rightarrow 0$ as $N\rightarrow \infty$, we can pick $N_0$ such that $8B_2 B_3 a_N \tau_N/3 \leq 1$ for any $N > N_0$. For any $\alpha > 0$, we can pick $M$ large enough so that $\frac{M^2}{8B_5 B_2 B_1 + M} \geq \alpha$ and $\exp\left(- \frac{M^2}{8B_5 B_2 B_1 + M} \ln N\right) < N^{\alpha}$. In particular, $M_{\alpha} = \frac{\alpha + \sqrt{\alpha^2 + 32B_5B_2B_1\alpha}}{2}$ will work. This means we have proved 
\begin{align*}
P\left(T_{N, 1}(w) > M_{\alpha} a_N\right) = O\left(N^{-\alpha}\right).
\end{align*}
We get the two-sided bound by applying the same argument twice for $T_{N, 1}(w)$ and $-T_{N, 1}(w)$. Moreover, because the derivation of the bound and the value of $M_\alpha$ doesn't depend on the specific point $w$, we have also proved our desired result
\begin{align*}
\sup_{w\in \RR^{d_W}}P\left(|T_{N, 1}(w)| > M_{\alpha} a_N\right) = O\left(N^{-\alpha}\right).
\end{align*}

\subsubsection*{Proof of claim (ii)}
We will use Propsition 2.3(c), a concentration inequality, from \cite{arcones1993}\footnote{There is a small modification compared to the original proposition. Since our statistic is not exactly a U-statistic as there are the iid $V_{ij}$ variables in our setup, we include this additional term in the statement of inequality. The proof of the inequality in our setup could follow the same steps of the original \cite{arcones1993} one. The reason this works is that the $V_{ij}$ terms are iid and won't affect the randomization inequality, decoupling inequality, and the hypercontractivity inequality used in the proof.} to prove the second claim.
\begin{quote}
Let $\{X_i, i\in \NN\}$ and $\{V_{i_1, \ldots, i_m}, (i_1, \ldots, i_m) \in I_m^\NN\}$ be independent random samples; $||f||_{\infty} \leq c$, $\EE \left] f(X_1, \ldots, X_m, V_{1, \ldots, m}) \right] = 0$, $\sigma^2 = \EE \left[ f^2(X_1, \ldots, X_m, V_{1, \ldots, m}) \right]$; $f$ is P-canonical, then there are constants $c_i$ depending only on $m$ such that for any $t>0$, 
\begin{multline*}
P\left(\left|N^{-m/2} \sum_{(i_1, \ldots, i_m) \in I_m^N} f(X_{i_1}, \ldots, X_{i_m}, V_{i_1, \ldots, i_m})\right| > t\right) \\
\leq c_1 \exp\left(-\frac{c_2 t^{2/m}}{\sigma^{2/m} + \left(c t^{1/m} N^{-1/2}\right)^{2/(m+1)}}\right).
\end{multline*}
\end{quote}
In order to apply the inequality, we first show that $\tau_N^{-1} h_N^{2d_X} \breve{Z}_{N,ij}$ is bounded. Decompose
\begin{align*}
\breve{Z}_{N, ij} & = \tilde{Z}_{N, ij} - \EE\left[\tilde{Z}_{N, ij} \Big |X_i, U_i\right] - \EE\left[\tilde{Z}_{N, ij} \Big |X_j, U_j\right] + \EE \tilde{Z}_{N, ij}.
\end{align*}
The last three terms on the right-hand side are bounded because $\tau_N^{-1}\EE \tilde{Z}_{N, ij} = O(1)$ and $\tau_N^{-1}\EE\left[\tilde{Z}_{N, ij} \Big |X_i, U_i\right] = O(h_N^{-d_X})$. Moreover,
$|\tau_N^{-1} h_N^{2d_X} \tilde{Z}_{N, ij}| = \frac{1}{2}
|\tau_N^{-1} Y_{ij}\cdot \indicator\left(|Y_{ij}| < \tau_N\right) K\left(\frac{w-W_{ij}}{h_N}\right)| + \frac{1}{2}
|\tau_N^{-1} Y_{ji}\cdot \indicator\left(|Y_{ji}| < \tau_N\right) K\left(\frac{w-W_{ji}}{h_N}\right)| \leq K_{\mathrm{max}}$. Hence, there exists constant $c>0$ s.t. $|\tau_N^{-1} h_N^{2d_X} \breve{Z}_{N,ij}| < c$. Applying the concentration inequality to $T_{N, 2}(w)$ then gives us
\begin{align*}
P\left(|T_{N, 2}(w)| > M a_N\right) & = P\left(\left|\frac{1}{\binom{N}{2}} \sum_{1 \leq i < j \leq N} \breve{Z}_{N, ij}\right| > M a_N\right)\\
& = P\left(\left|N^{-1} \sum_{1 \leq i < j \leq N} \tau_N^{-1} h_N^{2d_X} \breve{Z}_{N, ij} \right| > M \frac{N-1}{2}  h_N^{2d_X} a_N  \tau_N^{-1}\right)\\
& \leq c_1 \exp\left(-\frac{c_2 t}{\sigma + \left(c t^{1/2} N^{-1/2}\right)^{2/3}}\right)\\
& = c_1 \exp\left(-\frac{c_2 t}{\sigma \ln N+ \left(c t^{1/2} N^{-1/2}\right)^{2/3}\ln N} \cdot \ln N\right)\\
\end{align*}
where $t = M \frac{N-1}{2}  h_N^{2d_X} a_N  \tau_N^{-1}$ and $\sigma^2 = \var\left(\tau_N^{-1} h_N^{2d_X} \breve{Z}_{N, ij}\right)$. We will show that $\frac{c_2 t}{\sigma \ln N + \left(c t^{1/2} N^{-1/2}\right)^{2/3}\ln N} \rightarrow \infty$ as $N \rightarrow \infty$ by showing both $\frac{t}{\sigma \ln N} \rightarrow \infty$ and $\frac{t}{\left(c t^{1/2} N^{-1/2}\right)^{2/3}\ln N} \rightarrow \infty$ as $N \rightarrow \infty$. 

Beginning with the former claim:
\begin{align*}
\frac{t}{\sigma \ln N} & 
= \frac{M \frac{N-1}{2}  h_N^{2d_X} a_N  \tau_N^{-1}}{\tau_N^{-1} h_N^{2d_X} \var\left(\breve{Z}_{N, ij}\right)^{1/2} \ln N} 
= \frac{M (N-1) a_N}{2  \var\left(\breve{Z}_{N, ij}\right)^{1/2} \ln N} \geq \frac{M (N-1) a_N}{2 V_{N, 2}^{1/2} \ln N}\\
& \geq \frac{M N a_N}{4 \left(h_N^{-2d_X}B_4 K_{max} B_1\right)^{1/2} \ln N} = \frac{M}{4 \left(B_4 K_{max} B_1\right)^{1/2} } a_N \left(\frac{\ln N}{N h_N^{d_X}}\right)^{-1}\\
& =  \frac{M}{4 \left(B_4 K_{max} B_1\right)^{1/2} } a_N^{-1}\\
& \rightarrow \infty, \text{ as } N \rightarrow \infty. 
\end{align*}
The latter claim follows because:
\begin{align*}
\frac{t}{\left(c t^{1/2} N^{-1/2}\right)^{2/3}\ln N} & = \left(\frac{t^2 N}{c^2 (\ln N)^3}\right)^{1/3}
= \left(\frac{(M \frac{N-1}{2}  h_N^{2d_X} a_N  \tau_N^{-1})^2 N}{c^2 (\ln N)^3}\right)^{1/3}\\
& \geq \left(\frac{M^2}{16 c^2} N^3 (\ln N)^{-3} h_N^{4d_X} a_N^2  \tau_N^{-2}  \right)^{1/3}\\
& = \left(\frac{M^2}{16 c^2}\right)^{1/3}\left( N h_N^{\frac{3}{2}d_X} (\ln N)^{-1}\tau_N^{-1}\right)^{2/3}\\
& \rightarrow \infty, \text{ as } N \rightarrow \infty.
\end{align*}
The last line above is an implication of the condition $\tau_N \ll N h_N^{\frac{3}{2}d_X} /\ln N$. Combining these two limit results gives us $\frac{c_2 t}{\sigma \ln N + \left(c t^{1/2} N^{-1/2}\right)^{2/3}\ln N} \rightarrow \infty$ as $N \rightarrow \infty$. Notice the bound again doesn't depend on $w$ and the inequality still holds when we take the sup over $w \in \RR^{d_W}$ on the left-hand side. Hence for any $M>0$ and any $\alpha>0$, $\sup_{w \in \RR^{d_W}} P\left(|T_{N, 2}(w)| > M a_N\right) = O\left(N^{-\alpha}\right)$.

\subsubsection*{Proof of claim (iii)}
Direct evaluation yields
\begin{align*}
\left|\EE\left(\hat{\Phi}(w) - \tilde{\Phi}(w)\right)\right| 
& = \left|\EE \left[Y_{ij} \indicator\left(|Y_{ij}| > \tau_N\right) \frac{1}{h_N^{2d_X}} K\left(\frac{w-W_{ij}}{h}\right)\right]\right|\\
& \leq \EE \left[\left|Y_{ij}\right| \left|\tau_N^{-1} Y_{ij}\right|^{s-1}\indicator\left(|Y_{ij}| > \tau_N\right) \frac{1}{h_N^{2d_X}} \left|K\left(\frac{w-W_{ij}}{h_N}\right)\right|\right]\\
& \leq \tau_N^{-(s-1)} \EE \left[\left|Y_{ij}\right|^{s} \frac{1}{h_N^{2d_X}} \left|K\left(\frac{w-W_{ij}}{h}\right)\right|\right]\\
& = \tau_N^{-(s-1)} \int \EE \left[\left|Y_{12}\right|^{s} | (X_1, X_2)=(x_1, x_2)\right] \frac{1}{h_N^{2d_X}} \\
&\left|K\left(\frac{x - x_1}{h_N}, \frac{x - x_2}{h_N}\right)\right| f(x_1, x_2) \mathrm{d}x_1\mathrm{d}x_2\\
& = \tau_N^{-(s-1)} \int \EE \left[\left|Y_{12}\right|^{s} | (X_1, X_2)= (x-h_Ns_1, x'-h_Ns_2)\right] \\
& \times f(x - h_Ns_1, x'- h_Ns_2) \left|K\left(s_1, s_2\right)\right|  \mathrm{d}s_1\mathrm{d}s_2\\
& \leq \tau_N^{-(s-1)} B_{4,s} B_{1}.
\end{align*}
Since the last expression doesn't depend on $w$, we have $\sup_{w\in \RR^{d_W}}\left|\EE\left(\hat{\Phi}(w) - \tilde{\Phi}(w)\right)\right| = o(a_N)$.

\subsubsection*{Proof of claim (iv)}
First, we eliminate the sup by upper bounding the terms involving $K$ by $K_{\text{max}}$.
\begin{align*}
\sup_{w\in \RR^{d_W}}\left|\hat{\Phi}_N(w) - \tilde{\Phi}_N(w)\right| 
& = \sup_{w\in \RR^{d_W}}\left|\frac{1}{N(N-1)}\sum_{1\leq i \neq j\leq N} Y_{ij}\indicator\left(|Y_{ij}| > \tau_N\right) \frac{1}{h_N^{d_W}}K\left(\frac{w-W_{ij}}{h_N}\right)\right|\\
& \leq \frac{1}{N(N-1)}\sum_{1\leq i \neq j\leq N}\left|Y_{ij}\right| \indicator\left(|Y_{ij}| > \tau_N\right) \frac{1}{h_N^{2d_X}} \sup_{w\in \RR^{d_W}}\left|K\left(\frac{w-W_{ij}}{h_N}\right)\right|\\
& \leq K_{\text{max}} h_N^{-2d_X} \tau_N^{-(s-1)} \frac{1}{N(N-1)}\sum_{1\leq i \neq j\leq N}\left|Y_{ij}\right|^s.
\end{align*}
Then, taking expectation on both sides yields
\begin{align*}
\EE \left(\sup_{w\in \RR^{d_W}}\left|\hat{\Phi}_N(w) - \tilde{\Phi}_N(w)\right|\right) 
& \leq K_{\text{max}} h_N^{-2d_X} \tau_N^{-(s-1)}\EE\left(\left|Y_{ij}\right|^s\right) \leq K_{\text{max}} B_{6, s} h_N^{-2d_X} \tau_N^{-(s-1)} = o(a_N).
\end{align*}

\subsubsection*{Proof of claim (v)}
If all the $|Y_{ij}|, 1 \leq i \neq j \leq N$ are smaller than the truncation threshold $\tau_N$, then $\hat{\Phi}_N = \tilde{\Phi}_N$, 
\begin{align*}
P\left(\hat{\Phi}_N = \tilde{\Phi}_N\right) & \geq P\left(\max_{1\leq i < j\leq N}|Y_{ij}| \leq \tau_N \right).
\end{align*}
We now show that the RHS converges to $1$. Observe
\begin{align*}
\sum_{N=2}^\infty\sum_{i=1}^{N-1} [P\left(|Y_{iN}| > \tau_N \right) + P\left(|Y_{Ni}| > \tau_N \right)]
& \leq \sum_{N=2}^\infty\sum_{i=1}^{N-1} \left[\EE\left(|Y_{iN}|^s \tau_N^{-s} \right) + \EE\left(|Y_{Ni}|^s \tau_N^{-s} \right)\right]\\
& = \EE\left(|Y_{iN}|^s\right) \sum_{N=2}^\infty\sum_{i=1}^{N-1} 2N^{-2} \phi_N^{-1}\\
& \leq \EE\left(|Y_{iN}|^s\right) \sum_{N=2}^\infty \frac{2}{N (\ln\ln N)^2 \ln N}\\
& < \infty,
\end{align*}
The Borel-Cantelli lemma implies $P(A_{ij}, i\neq j, i.o.) = 0$ where the set $A_{ij} = \{\omega: Y_{ij}(\omega) > \tau_{\max\{i, j\}}\}$. This means, except for a null set $\cN$, for any $\omega \in \cN^c$, there exists a $N(\omega)$ s.t. for all $N \geq N(\omega)$, $Y_{iN}(\omega) \leq \tau_N$. Since $\tau_N \uparrow \infty$ as $N \rightarrow \infty$, we can take $N^*(\omega) \geq N(\omega)$ such that $\tau_{N^*(w)} > \max_{i,j \leq N(\omega)} |Y_{ij}(\omega)|$. Then for any $N \geq N^*(\omega)$, we have $\max_{1 \leq i < j \leq N}|Y_{ij}(\omega)| \leq \tau_N$ and hence $\hat{\Phi}_N = \tilde{\Phi}_N$. Define the set $E_N := \{\omega: N^*(\omega) \leq N\} \subset \{\omega: \hat{\Phi}_N = \tilde{\Phi}_N\}$. Since $E_N \uparrow \cN^c$ and $P(\cN^c) = 1$, we have $P(\hat{\Phi}_N = \tilde{\Phi}_N) \geq P(E_N)\rightarrow 1$ as $N \rightarrow \infty$.

\subsection*{Proof of Theorem \ref{thm:UniformConvergenceRegression}}
The proof follows the general approach used in \cite{hansen2008}.
Denote 
$\hat{f}_{W, N}(w) = \frac{1}{N(N-1)}\sum_{1\leq i \neq j \leq N} K_{ij, N}(w)$. We can write
\begin{align*}
\hat{g}_N(w)
& = \frac{\hat{\Psi}_N(w)}{\hat{f}_{W, N}(w)}.
\end{align*}
We examine the numerator and denominator separately. An application of Theorem \ref{thm:UniformConvergence} yields
\begin{align*}
\sup_{||w||\leq C_N} |\hat{\Psi}_N(w) - \EE\hat{\Psi}_N(w)| 
& = O_p(a_N)\\
\sup_{||w||\leq C_N} |\hat{f}_{W, N}(w) - \EE\hat{f}_{W, N}(w)| 
& = O_p(a_N).
\end{align*}
Standard bias calculations give
\begin{align*}
\sup_{||w||\leq C_N} |\EE\hat{\Psi}_N(w) - \Psi(w)| 
& = O(h_N^{\beta})\\
\sup_{||w||\leq C_N} |\EE\hat{f}_{W, N}(w) - f_W(w)| 
& = O(h_N^{\beta}).
\end{align*}
Combining these results we get
\begin{align*}
\sup_{||w||\leq C_N} |\hat{\Psi}_N(w) - \Psi(w)| 
& = O_p(a_N) + O(h_N^{\beta}) = O(a_N^*)\\
\sup_{||w||\leq C_N} |\hat{f}_{W, N}(w) - f_W(w)| 
& = O_p(a_N) + O(h_N^{\beta}) = O(a_N^*).
\end{align*}

Uniformly over $||w|| \leq C_N$ we have
\begin{align*}
\frac{\hat{\Psi}_N(w)}{\hat{f}_{W, N}(w)}
& = \frac{\hat{\Psi}_N(w)/f_W(w)}{\hat{f}_{W, N}(w)/f_W(w)}
= \frac{g(w) + (\hat{\Psi}_N(w) - \Psi(w))/f_W(w)}{1 + (\hat{f}_{W, N}(w) - f_W(w))/f_W(w)}
= \frac{g(w) + O_p(\delta_N^{-1}a_N^*)}{1 + O_p(\delta_N^{-1}a_N^*)}\\
& = g(w) + O_p(\delta_N^{-1}a_N^*)
\end{align*}
as claimed. The optimal rate is obtained by setting $h_N \asymp \left(\frac{\ln N}{N}\right)^{\frac{1}{2\beta+d_X}}$.

\end{document}